\documentclass{amsart}

\usepackage{amsmath,amsfonts,amssymb}
\usepackage{amsthm,amsfonts,amssymb}
\usepackage{graphics, verbatim}
\usepackage{pifont}
\usepackage[usenames,dvips]{color}
\usepackage{amsrefs}
\def\({\left(}
\def\){\right)}
\def\be{\begin{equation}}
\def\en{\end{equation}}

\let\nc\newcommand
\nc{\al}{\alpha}

\nc{\BB}{{A_{p,q}}}
\nc{\bean}{\begin{eqnarray}} \nc{\eean}{\end{eqnarray}}
\nc{\bea}{\begin{eqnarray*}} \nc{\eea}{\end{eqnarray*}}

\nc{\bs}{\boldsymbol}

\def\Br{\mathbb R}

 \nc{\ii}{{ i}}
 \nc{\pp}{{ p}}

\newcommand{\Z}{\mathbb{Z}}

\theoremstyle{plain}
\begingroup 
\newtheorem{thm}{Theorem}[section]   
\newtheorem*{thm*}{Theorem}          
\newtheorem*{cor*}{Corollary}        
\newtheorem{cor}[thm]{Corollary}     

\newtheorem{lem}[thm]{Lemma}         
\newtheorem{prop}[thm]{Proposition}  
\endgroup

\theoremstyle{definition}

\newtheorem*{rem*}{Remark}
\newtheorem*{ack*}{Acknowledgment}

\theoremstyle{remark}
\newtheorem{ex}[thm]{Example}        

\theoremstyle{definition}

\numberwithin{equation}{section}

\begin{document}

\title[Path counts and Stirling numbers]
{Path count asymptotics and Stirling numbers}

\author{K. Petersen}
\address{Department of Mathematics, University of North Carolina
at Chapel Hill\\ Chapel Hill, NC 27599-3250, USA} \email{petersen@math.unc.edu}

\author{A. Varchenko}
\address{Department of Mathematics, University of North Carolina
at Chapel Hill\\ Chapel Hill, NC 27599-3250, USA} \email{anv@math.unc.edu}
\thanks{The research of
the second author was supported in part by NSF grant DMS-0555327}

\date{\today}


\begin{abstract}
{
We obtain formulas for the growth
rate of the numbers of certain paths in a multi-dimensional analogue of the Eulerian graph.
Corollaries are
new identities} relating Stirling numbers of the first and second kinds.
\end{abstract}

\maketitle

\section{Introduction}
\label{sec-intro}

{
A graded, infinite, directed graph defines a dynamical system called a
Bratteli-Vershik or adic system (see \cite{PetersenVarchenko2009} and the references given there for background on adic systems). One of the main problems concerning these systems is to determine
their ergodic invariant measures, and the basic tool for attacking this problem
is the asymptotics of numbers of paths in the graph with certain
properties.}

{
 A very interesting Bratteli-Vershik dynamical system
is associated to the Eulerian two-dimensional graph with vertices at the points of
 $\Z^2_{\geq 0}$. The paths in the Eulerian graph may be counted
in many different combinatorial ways, but}
{
all attempts to extend those counts
to the multi-dimensional analogues}
{
 of the Eulerian graph have failed
up to this moment;
techniques applied to the
   two-dimensional Eulerian graph, such as the recurrence of random
   walks in two dimensions in \cite{BKPS}, the coding of paths by
   permutations in \cite{FP}, or the left-to-right monotonicity in
   \cite{Gnedin-O}, are not available in dimensions higher than two.}

In this paper we make a first step in the multi-dimensional direction.
We consider a multi-dimensional analogue of the Eulerian graph and determine the
asymptotics of the numbers of paths with certain properties.
We plan to apply these asymptotics to studying the ergodic
measures on these and related systems in future work.
As a rather unexpected
{corollary}
 of our asymptotics we get a family of
{
 new
identities} between Stirling numbers of the first and second kinds, which
are classical objects,
arise in a variety of combinatorics problems, and
 for which many general relations have long been known (see for example \cites{Sloane2010,SloanePlouffe1995}).

 {Now we  formulate the results
  of this paper}.
 Denote by ${} e_1=(1,0,\dots,0),$ \dots , ${}
e_{n+1}=(0,\dots,0,1) \in \Z^{n+1}$ the ordinary basis vectors in
$\Br^{n+1}$.
   Let $c \in \Z_{\geq 0}^{n+1}$.
 We define the graph $E_c$ to be the graph with vertices $\ii\in \Z_{\geq 0}^{n+1}$ and the
following directed multiple edges.
 Each vertex $\ii$ is connected just to $n+1$
vertices:
 to vertex $\ii  +e_j$ by $c_j+i_{n+1}$ edges if $j<n+1$ and to $\ii +e_{n+1}$  by
 $c_{n+1}+i_1+i_2+\dots + i_{n}$ edges.

Such a graph models a random walk on the graph $G_{n+1}$ consisting of
$n+1$ loops based at a single
 vertex, with a special kind of opposite reinforcement:
 when the walker takes a step in the direction of $e_{n+1}$,
 the numbers of edges in all the other directions ($e_j$
 for $j=1, \dots , n$) are incremented by 1; but when the walker
 takes a step in the direction of $e_j$ for some $j=1, \dots , n$,
 {\em only} the number of edges in the direction of $e_{n+1}$ is
 incremented by $1$ (cf. \cite{FP2}).
 When $n=1$ and $c=(1,1)$, $E_c$ is the
 {standard}
two-dimensional Eulerian
 graph.
 The corresponding two-dimensional Eulerian adic system
 was introduced in \cite{BKPS} and studied further in \cite{FP} and
 \cite{PetersenVarchenko2009}; see also \cites{CarlitzScoville1974,
 Flajolet2006, Gnedin-O}.

 For each $i \in \Z_{\geq 0}^{n+1}$ denote by $A_c(i)$
 the number of paths in the graph $E_c$ from $0$ to $i$.
 The standard Eulerian numbers
 (see \cites{Comtet, PetersenVarchenko2009}) are the numbers
 $A(i_1,i_2)=A_{(0,0)}(i_1,i_2)$
{ with $n=1$.
In this paper, we give a formula for the leading term of the growth rate of
$A_c(i_1,\dots,i_n,i_{n+1})$ as $i_{n+1}$ tends to infinity.}
We show that if $m=i_1+\dots + i_n$, then
$$
A_c(i_1,\dots,i_n,i_{n+1}) \sim
B_c(i_1,\dots,i_n) (c_{n+1}+ m)^{i_{n+1}}
$$
as $i_{n+1} \to \infty$ and \label{introduction prop-1 and i}
 \be \label{B(i)}
 B_{c_1,\dots, c_{n+1}}(i_1,\dots,i_n)\ =\
 \frac {( c_{n+1}+m)^{m}}{m!} \, \sum_{f\in F} \,
 \prod_{j=1}^{m}\ (c_{f(j)} + c_{n+1} + j-1)\ ,
 \en
where
$F$ denotes the set of
all maps $f:\{1,\dots,m\}\to \{1,\dots,n\}$ such that $|f^{-1}(j)| =
i_j$ for   $j=1,\dots,n$.

As a corollary we obtain {new identities} relating
 Stirling numbers of the first and second
 kinds;  for
$1\leq k \leq n$,\ $0 \leq r \leq k$, we get
  \be
\label{introductioncoefs}
  \begin{gathered}
  {r+n-k-1 \choose r}
  s_1(n,r+n-k)=\\
  \sum_{m=0}^k {m+n-k \choose m+1}\sum_{i=0}^r
  {i+n-k+m-1 \choose i}\frac{(-1)^{m+r-i}}{n^{r-i+1}}(r-i+1)!\ \, \times \\
  s_2(m+1,r-i+1) s_1(n,i+n-k+m) ,
  \end{gathered}
  \en
where $s_1, s_2$ are Stirling numbers of the first and second kinds:
see Theorem \ref{thm S s} and Corollary
\ref{cor-stirlingidentities}.

\begin{ack*}
We thank R. Stanley for telling us about theorem of Frobenius
(Theorem \ref{stanley}) and the referee for suggestions leading to the improvement of the exposition.
\end{ack*}

\section{Path counts and limits}
\subsection{Finding the leading term}\label{subsec-multidim}

 The numbers $A_c(i)$ of paths in the graph $E_c$ from 0 to $i$
 satisfy the recurrence relation
 \be
 A_{c}(\ii) = (c_{n+1}+i_1+\dots +i_n) A_c(i-e_{n+1}) +
 \sum_{j=1}^{n}(c_j + i_{n+1}) A_c(\ii-e_j) .
 \en

We are interested in the limit
 \be\label{eq-Bc}
 B_c(i_1,\dots,i_{n})=\lim_{i_{n+1} \to \infty}
 \frac{A_c(i_1,\dots,i_{n},i_{n+1})
}{(c_{n+1}+n)^{i_{n+1}}} .
 \en


\begin{prop}\label{prop-Bform}
Let $u=(1,\dots, 1)\in \mathbb Z_{\geq 0}^n.$ Then
  \be
\label{B}
 B_c(u)  =
\sum_{\xi\in S_{n}} \sum_{j_1,\dots,j_{n}=0}^\infty\,
\prod_{k=1}^{n}\ (c_{\xi(k)} +j_1+\dots+j_k) \left(\frac
{c_{n+1}+k-1}{c_{n+1}+{n}}\right)^{j_k} .
 \en
   \end{prop}
  \begin{proof}
  To form a path from ${} 0$ to $\ii = (1, \dots, 1,i_{n+1})$ in
  $E_c$, we need to take one step in the direction of each ${}
  e_k$,
  $k=1,\dots, n$,
  and $i_{n+1}$ steps in the direction of ${} e_{n+1}$. Thus each such
  path will correspond to a permutation $\xi$ of $1,2,\dots, n$
  which will determine when the step in the direction of each ${} e_k, k=1,\dots,n,$
  is taken, together with a partition
  $i_{n+1}=j_1+\dots+j_n+j_{n+1}$ of $i_{n+1}$ into nonnegative
  integers, which will determine how many consecutive steps are
  taken in the ${} e_{n+1}$ direction between steps in other directions.
  Each time that $j$ steps are taken in the ${} e_{n+1}$ direction, all the
  $c_k, 1 \leq k \leq n,$ are incremented by 1; and each time that a
  step is taken in a direction other than ${} e_{n+1}$, $c_{n+1}$
  is incremented by 1.
  Therefore,
  \be
  \begin{aligned}
  &A_c (1, \dots, 1,i_{n+1})
\\
   &=  \sum_{\xi \in
  S_n}\sum_{j_1+\dots + j_{n+1}=i_{n+1}}
  c_{n+1}^{j_1}(c_{\xi(1)}+j_1)(c_{n+1}+1)^{j_2}(c_{\xi(2)}+j_1+j_2) \cdots
\\
  &\hphantom{aaaaaaa}\cdots
  (c_{\xi(n)}+j_1+\dots +j_n)(c_{n+1}+n)^{j_{n+1}}
\\
  &= \sum_{\xi \in  S_n}\sum_{j_1+\dots + j_{n+1}=i_{n+1}}
  \left[ \prod_{k=1}^n(c_{n+1}+k-1)^{j_k}(c_{\xi(k)}+j_1+\dots
  +j_k)\right] (c_{n+1}+n)^{j_{n+1}}
\\
  &= (c_{n+1}+n)^{i_{n+1}}
\!\! \sum_{\xi \in  S_n}\sum_{j_1+\dots + j_{n+1}=i_{n+1}}
  \prod_{k=1}^n \left[ \frac{(c_{n+1}+k-1)}{c_{n+1}+n} \right]^{j_k}
  \!\!
(c_{\xi(k)}+j_1+\dots +j_k) \\
  &= (c_{n+1}+n)^{i_{n+1}}\sum_{\xi \in
  S_n}\sum_{\genfrac{}{}{0 pt}{}{j_s \geq 0}{j_1+\dots + j_{n}=0}}^{i_{n+1}}
  \prod_{k=1}^n \left[ \frac{(c_{n+1}+k-1)}{c_{n+1}+n} \right]^{j_k}
  (c_{\xi(k)}+j_1+\dots +j_k).
  \end{aligned}
  \en

This reasoning gives Formula (\ref{B}). Clearly, $B_c(u)$
 is a convergent sum.
 \end{proof}

 Similarly, existence of the limit in Formula
 (\ref{eq-Bc}) for any $(i_1,\dots,i_n)$ can be shown as follows.
 Let $m$, $F$, be the same as in Section \ref{sec-intro}. Then
 \be
 \begin{aligned}
 &A_c(i_1,\dots,i_n,i_{n+1})=\\
   &=  \sum_{f \in F}\sum_{j_1+\dots + j_{m+1}=i_{n+1}}
  c_{n+1}^{j_1}(c_{f(1)}+j_1)(c_{n+1}+1)^{j_2}(c_{f(2)}+j_1+j_2) \cdots
\\
  &\hphantom{aaaaaaa}\cdots
  (c_{f(m)}+j_1+\dots +j_m)(c_{n+1}+m)^{j_{m+1}}
\\
   &= (c_{n+1}+m)^{i_{n+1}}\sum_{f \in F}
\sum_{\genfrac{}{}{0 pt}{}{j_s \geq 0}{j_1+\dots + j_{m}=0}}^{i_{n+1}}
  \prod_{k=1}^m \left[ \frac{(c_{n+1}+k-1)}{c_{n+1}+m} \right]^{j_k}
  (c_{f(k)}+j_1+\dots +j_k).
  \end{aligned}
  \en
Dividing by $(c_{n+1}+m)^{i_{n+1}}$ and letting $i_{n+1} \to \infty$
again leaves a sum of convergent series. Thus we get
\be
\label{B i}
 B_c(i_1,\dots,i_n)  =
\sum_{f\in F}
\sum_{j_1,\dots,j_{n}=0}^\infty\,
\prod_{k=1}^{n}\ (c_{f(k)} +j_1+\dots+j_k) \left(\frac
{c_{n+1}+k-1}{c_{n+1}+{n}}\right)^{j_k} .
 \en


We shall give two formulas for $B_c(u)$. In this section we evaluate
$B_c(u)$ using a {\em two-dimensional} formula from
 \cite{PetersenVarchenko2009}.
In Section
 \ref{subsec-derivsofseries} another formula is obtained by
using derivatives of geometric series. Comparing the two formulas we
get new  identities involving Stirling numbers of the first and
second kinds
 (Corollary \ref{cor-stirlingidentities}).

\begin{thm}
 \label{conj} We have
 \be \label{C}
 B_c(u)\ =\ \frac
{(c_{n+1}+n)^{n}}{n!} \, \sum_{\xi\in S_{n}} \, \prod_{k=1}^{n}\
(c_{\xi(k)} + c_{n+1} + k-1)\ . \en
\end{thm}

\subsection{Proof of Theorem \ref{conj}} \label{subsec-pfofconj}

Denote by $\sigma_i({} c),\ i=1,\dots,n$, the elementary symmetric
functions in ${} c$, \be \sigma_i({} c) = \sum_{1\leq j_1<\dots
<j_i\leq n}c_{j_1}\dots c_{j_i} \ . \en Denote $\sigma_0({} c)=1$.
Clearly, the function $B_c(u)$ is a symmetric polynomial in ${} c$
of
 the form
 \be \label{Rn}
 \begin{gathered}
 B_c(u) = \sum_{i=0}^n\ i!(n-i)!\ \sigma_{i}({} c)
 \al_{i,n}(c_{n+1}) =n!\,\sigma_n({} c) \al_{n,n}(c_{n+1}) +\\
   (n-1)!\,\sigma_{n-1}({}
c) \al_{n-1,n}(c_{n+1})\ +\ \dots\ +\ n! \,\sigma_0({} c)
\al_{0,n}(c_{n+1})\ ,
\end{gathered}
 \en
 where   $\al_{i,n}(c_{n+1}),\
i=0,\dots,n$, are suitable functions in $c_{n+1}$. For example,
 \be
 \begin{aligned}
\al_{0,n}(c_{n+1}) &= \sum_{j_1,\dots,j_{n}=0}^\infty\,
\prod_{k=1}^{n}\ (j_1+\dots+j_k) \left(\frac
{c_{n+1}+k-1}{c_{n+1}+{n}}\right)^{j_k} \ ,\\
 \al_{n,n}(c_{n+1})  &=\ \sum_{j_1,\dots,j_{n}=0}^\infty\,
\prod_{k=1}^{n}\ \left(\frac
 {c_{n+1}+k-1}{c_{n+1}+{n}}\right)^{j_k},
 \end{aligned}
 \en
 and so on.

\medskip
 \begin{ex}
 For $n=3$,
 \be
 \begin{aligned}
 B_{(c_1,c_2,c_3)}(u) &=\
\sum_{j_1,j_2=0}^\infty
[(c_1+j_1)(c_2+j_1+j_2)+(c_2+j_1)(c_1+j_1+j_2)]\ \times
\\
&\phantom{aaa} \times  \left(\frac{c_3}{c_3+2}\right)^{j_1}
\left(\frac{c_3+1}{c_3+2}\right)^{j_2}
\\
 &= \ 2!\, c_1c_2 \,\al_{2,2}(c_3) \ +\ (c_1+c_2)\,\al_{1,2}(c_3)\ +\  2! \al_{0,2}(c_3) ,
 \end{aligned}
 \en
 where
 \be  \label{al}
 \begin{aligned}
  \al_{2,2}(c_3) &=\sum_{j_1,j_2=0}^\infty \left(\frac{c_3}{c_3+2}\right)^{j_1}
 \left(\frac{c_3+1}{c_3+2}\right)^{j_2} ,\\
  \al_{1,2}(c_3) &= \sum_{j_1,j_2=0}^\infty
(2j_1+j_2)\left(\frac{c_3}{c_3+2}\right)^{j_1}
 \left(\frac{c_3+1}{c_3+2}\right)^{j_2} ,\\
  \al_{0,2}(c_3) &= \sum_{j_1,j_2=0}^\infty
j_1(j_1+j_2)\,\left(\frac{c_3}{c_3+2}\right)^{j_1}
\left(\frac{c_3+1}{c_3+2}\right)^{j_2} .
\end{aligned}
 \en
 \end{ex}

Our goal is
 to calculate all the coefficients $\al_{i,n}(c_{n+1})$ and show that the result gives
Formula (\ref{C}).
 To do so let us calculate
 the number
 \be \label{eq_secondproduct}
B_{c_1,c_2}(n) \ =  \ \sum_{j_1,\dots,j_{n}=0}^\infty\,
\prod_{k=1}^{n} \ (c_1 + j_1+\dots + j_k) \left(\frac
{c_{2}+k-1}{c_{2}+{n}}\right)^{j_k} .
 \en
  The function
 $B_{c_1,c_2}(n)$ has the form:
 \be
 B_{c_1,c_2}(n) \ =\ c_1^n
 \al_{n,n}(c_2) + c_1^{n-1} \al_{n-1,n}(c_2) + \dots +
 \al_{0,n}(c_2) ,
 \en
 where the functions $\al_{i,n}$, $i=0,\dots,n$,
are the {\em same} functions of one variable as in
 (\ref{Rn}).
  (To see this, note the identical form of the two
 products in formulas (\ref{B}) and (\ref{eq_secondproduct}), and note that
the sum over the symmetric group is absorbed into the coefficients
that do not involve $c_{n+1}$.)

\medskip
 \begin{ex}
  We have
 \be
 \begin{aligned}
  B_{c_1,c_2}(2) &=
 \sum_{j_1,j_2=0}^\infty (c_1+j_1)(c_1+j_1+j_2)\left(\frac{c_2}{c_2+2}\right)^{j_1}
 \left(\frac{c_2+1}{c_2+2}\right)^{j_2}
\\
 &= c_1^2 \,\al_{2,2}(c_2) \ +\ c_1\,\al_{1,2}(c_2)\ +\  \al_{0,2}(c_2) ,
 \end{aligned}
 \en
 where the functions $\al_{0,2},\al_{1,2},\al_{2,2}$ are the same
functions of one variable as in (\ref{al}).
 \end{ex}

To calculate the functions $\al_{i,n}$ and prove Theorem
\ref{conj}
it is enough to calculate $B_{c_1,c_2}(n)$.

\begin{lem}
\label{lem main}
We have \be \label{1} B_{c_1,c_2}(n)\ =\ \frac
{(c_{2}+n)^{n}}{n!} \, \prod_{k=1}^{n}\ (c_{1} + c_{2} + k-1)\ . \en
\end{lem}

\begin{proof}
Let us use formula (1.1) of \cite{PetersenVarchenko2009} to evaluate
 \be \label{eq-1}
B_{c_1,c_2}(n)=\lim_{i\to\infty}\frac
 {A_{c_1,c_2}(n,i)}{(c_2+n)^i}.
 \en
 Formula (1.1) for
$A_{c_1,c_2}(n,i)$ has $n+1$ summands and the leading summand is
 \be
{c_1+c_2+n-1\choose n}\,(c_2+n)^{n+i}.
 \en
  This term after division by
$(c_2+n)^i$ gives the right hand side of (\ref{eq-1}). The lemma is
proved.
\end{proof}

Lemma \ref{lem main}  gives Theorem \ref{conj}.

\subsection{Formula for
 $B_c(\ii)$}\label{subsec-generalB}
 Formula (\ref{C}) for
$ B_{c_1,\dots,c_{n+1}}(1,\dots,1)$  easily gives a formula for
 $ B_{c_1,\dots,c_{n+1}}(i_1,\dots,i_n)$
for any $(i_1,\dots,i_n)$ .




\begin{prop}
\label{prop-1 and i}
 Let $m$, $F$, be the same as in Section \ref{sec-intro}. Then
 \be \label{B(i)}
 B_{c_1,\dots, c_{n+1}}(i_1,\dots,i_n)\ =\
 \frac {( c_{n+1}+m)^{m}}{m!} \, \sum_{f\in F} \,
 \prod_{j=1}^{m}\ (c_{f(j)} + c_{n+1} + j-1)\ .
 \en
\end{prop}
\bigskip

\section{Stirling numbers and related polynomials}\label{sec-StirlingNumbers}

Theorem \ref{conj} has as corollaries interesting
identities involving Stirling numbers of the first and second kinds.

 \subsection{Theorem \ref{conj} in new notation}
\label{subsec-derivsofseries}
 Denote
\be u_a=\frac{q-a}q\ , \qquad D_a = u_a\frac\partial{\partial u_a}\
. \notag \en Then Theorem \ref{conj} says that for arbitrary numbers
$c_1,\dots,c_n$
 we have
\be \label{new n}
 \left( \sum_{\sigma\in S_n} \prod_{a=1}^n \,(c_{\sigma(a)} +
D_a+D_{a+1}+\dots+D_n)\! \right) \!\prod_{b=1}^n\frac 1{1-u_b}  =
\frac {q^n}{n!}\!\sum_{\sigma\in S_n} \prod_{a=1}^n (c_{\sigma(a)}
 + q-a).
 \en

\begin{ex}
For $n=2$, Formula  (\ref{new n}) says
 \be
\begin{aligned}
& [(c_1+D_2)(c_2+D_2+D_1)+(c_2+D_2)(c_1+D_2+D_1)]\frac 1{1-u_1}\frac
1{1-u_2}\ = \
\\
& \phantom{aaaaaaaaaaaaaaaaaa} \frac {q^2}2 \ [(c_1 + q-1)(c_2 +
 q-2)+(c_2 + q-1)(c_1 + q-2)]\ .
 \end{aligned}
\notag \en
 This statement is equivalent to
 the identities
 \bea
 \frac 1{1-u_1}\frac 1{1-u_2}\ &=& \ \frac {q^2}2\
,
\\
(2D_2+D_1)\frac 1{1-u_1}\frac 1{1-u_2}\ &=& \ \frac {q^2}2\
(q-2+q-1)\ ,
\\
D_2(D_2+D_1)\frac 1{1-u_1}\frac 1{1-u_2}\ &=& \ \frac {q^2}2\
 (q-1)(q-2)\ .
 \eea
\end{ex}


 \subsection{Formula for $(u\partial_u)^k(1-u)^{-1}$}\label{subsec-derivforms}

We make some observations that are useful for handling expressions
such as the left side of Formula (\ref{new n}).

 Let $q$ be a variable, $a$ a positive integer. Set
 \be
\label{eq-diff}
  u =\frac{q-a}q \ , \qquad \partial_u = \frac{\partial}{\partial u}\ .
 \en
Let $A(i,j)=A_{0,0}(i,j)$ be the standard Eulerian numbers as in
Section \ref{sec-intro}.

\begin{thm}
\label{thm 1} For any positive integer $k$
we have \be \label{form 1}
(u\partial_u)^{k}(1-u)^{-1}\ = \ a^{-k-1}\!\!\sum_{i+j=k-1}
A(i,j)\,q^{i+1}(q-a)^{j+1}\ . \en
\end{thm}

\begin{ex}
 \be
 \begin{aligned}
 u\partial_u(1-u)^{-1}\ &= \
a^{-2}\,q(q-a)\ ,\\
 (u\partial_u)^2(1-u)^{-1}\ &= \ a^{-3}\,(q^2(q-a) + q(q-a)^2)\ ,\\
 (u\partial_u)^3(1-u)^{-1}\ &= \ a^{-4}\,(q^3(q-a) + 4 q^2(q-a)^2 +
 q(q-a)^3)\ .
\end{aligned}
\en
\end{ex}

Theorem \ref{thm 1} follows by induction from the following lemma.

\begin{lem}
\label{lem der} We have \be
 u\frac {\partial}{\partial u}\  =\ \frac
{q(q-a)}a \frac {\partial}{\partial q} \ . \notag \en
\end{lem}

\noindent
 {\em Proof of Theorem \ref{thm 1}.}
 \be
 \begin{gathered}
 \frac {q(q-a)}a \frac {\partial}{\partial q}
 a^{-k-1}\sum_{i+j=k-1}
A(i,j)\,q^{i+1}(q-a)^{j+1}
\\
 = a^{-k-2}\sum_{i+j=k-1}A(i,j)\,(i+1)\, q^{i+1}(q-a)^{j+2}
 \\
+\, a^{-k-2}\sum_{i+j=k-1}A(i,i)\,(j+1)\, q^{i+2}(q-a)^{j+1}
\\
 =a^{-k-2}\sum_{r+s=k}A(r,s-1)\,(r+1)\, q^{r+1}(q-a)^{s+1}
\\
 + \,a^{-k-2}\sum_{r+s=k}A(r-1,s)\,(s+1)\, q^{r+1}(q-a)^{s+1}
\\
 = \,a^{-k-2}\sum_{r+s=k}A(r,s)\, q^{r+1}(q-a)^{s+1}\ .
 \end{gathered}
\notag
 \en
\qed

\subsection{$\Gamma$-polynomials}
\label{subsec-gammas}

We shall reformulate Theorem \ref{thm 1} in terms of a slight
variation of Eulerian polynomials (see \cite{Comtet}*{p. 244}),
which we call $\Gamma$-polynomials.

  Let $s_2(k,m)$ denote
Stirling numbers of the second kind,
{which count the number of partitions of a $k$-element set into
$m$ nonempty subsets, so that }
\bea &s_2(1,1)=1, &
\\
&s_2(2,1)=1, \qquad s_2(2,2)=1, &
\\
 &s_2(3,1)=1, \qquad s_2(3,2)=3, \qquad s_2(3,3)=1.
&
 \eea
 We have
 \be
 s_2(k,m) = s_2(k-1,m-1) + m s_2(k-1,m)\ .
 \notag
\en
 We will be interested in the polynomials
 \be
 \begin{aligned}
 \Gamma_k(q,n)\ & =\
\sum_{i=1}^k\,(-1)^{k-i} \,s_2(k,i) \,i!\,q^{i-1}n^{k-i}\ .
 \end{aligned}
\notag
 \en

 \begin{ex}
 \be
\Gamma_1(q,n) = 1 , \qquad \Gamma_2(q,n) = 2q-n , \qquad
\Gamma_3(q,n) = 6q^2-6qn +n^2 . \notag \en
 \end{ex}

\begin{thm} [Frobenius, \cite{Comtet}*{p. 244}]
\label{stanley} We have \be \Gamma_k(q,n)\ = \ \sum_{i+j=k-1}
A(i,j)\,q^{i}(q-n)^{j}\ . \notag \en
\end{thm}

\subsection{$\Delta$-polynomials} \label{subsec-taus}
We will also need $\Delta$-polynomials.

 Looking at the right side of Formula (\ref{new n}) suggests
 investigating related polynomials and
 Stirling numbers of the first kind.

 Stirling numbers of the first kind are the coefficients in the
expansion \be
    (q)_{n} = \sum_{k=0}^n s_1(n,k) q^k,
\notag \en where $(q)_n$ is the falling factorial
 \be
    (q)_{n}=q(q-1)(q-2)\cdots(q-n+1) .
 \notag
\en
 We have
 \be
s_1(n,k) = \sum_{1\leq i_1<\dots<i_{n-k}\leq n-1}
 (-1)^{n-k} \,i_1i_2\cdots i_{n-k} \ .
\notag
 \en
 For
natural numbers $1\leq k < n$, define a polynomial
$$
\Delta_{n,k}(q) = \sum_{1\leq i_1<\dots<i_{n-k}\leq n-1}
(q-i_1)(q-i_2)\cdots(q-i_{n-k}) .
$$
Define $\Delta_{n,n}(q) = 1$. We have $\Delta_{n,1}(q)\, =\,
(q-1)\dots (q-n+1)$\  and \be \Delta_{n,k}(q)\ = \
\sum_{m=0}^{n-k}\, {k-1+m\choose m}\,s_1(n,k+m)\,q^{m}\ . \notag \en

\begin{ex}
 \be
 \Delta_{2,2}(q) = 1, \qquad        \Delta_{2,1}(q) = q-1 ,
\notag \en \be \Delta_{3,3}(q)=1, \qquad       \Delta_{3,2}(q) =
2q-3, \qquad \Delta_{3,1}(q)=q^2-3q+2 . \notag \en
\end{ex}

\subsection{Relation between $\Gamma$ and $\Delta$ polynomials}
 \label{subsec-polyform}

 \begin{ex}
 We have the following identity:
 \be
\notag  \left(
\begin{array}{clcr} {}^{2^{-1}{2\choose 1} {\Gamma_{1}(q,2)}}
\phantom{aaaa} & \phantom{aaaa} {}^{0}
\\
{}^{2^{-2}{2\choose 2}\Gamma_{2}(q,2)}\phantom{aaaa} &
{}^{2^{-1}{1\choose 1}\Gamma_{1}(q,2)}
\end{array} \right)
\left( \begin{array}{clcr} \Delta_{2,2}(q)
\\
\Delta_{2,1}(q)
\end{array} \right)
= \left( \begin{array}{clcr} \Delta_{2,2}(q)
\\
\Delta_{2,1}(q)
\end{array} \right) .
\en In other words, we have
  \be
\notag
  {2\choose
1}\frac{\Gamma_{1}(q,2)\, \Delta_{2,2}(q)}2 = \Delta_{2,2}(q)\ , \en
\be {2\choose 2}\frac{\Gamma_{2}(q,2)\, \Delta_{2,2}(q)}{2^2}
 + {1\choose 1}\frac{\Gamma_{1}(q,2)\, \Delta_{2,1}(q)}{2}
 = \Delta_{2,1}(q)\ .
\notag\en
 \end{ex}

 \begin{ex}
 We have the following identity
\be \notag
 \left(
\begin{array}{clcr} {}^{3^{-1}{3\choose 1}\Gamma_{1}(q,3)}
\phantom{aaaa} & \phantom{aaa} {}^{0} & \phantom{aaaa} {}^{0}
\\
{}^{3^{-2}{3\choose 2}\Gamma_{2}(q,3)}\phantom{aaaa} &
{}^{3^{-1}{2\choose 1}\Gamma_{1}(q,3)} & \phantom{aaaa}{}^{0}
\\
{}^{3^{-3}{3\choose 3}\Gamma_{3}(q,3)} \phantom{aaaa} &
{}^{3^{-2}{2\choose 2}\Gamma_{2}(q,3)} & \phantom{aaaa}
{}^{3^{-1}{1\choose 1}\Gamma_{1}(q,3)}
\end{array} \right)
\left( \begin{array}{clcr} \Delta_{3,3}(q)
\\
\Delta_{3,2}(q)
\\
\Delta_{3,1}(q)
\end{array} \right)
= \left( \begin{array}{clcr} \Delta_{3,3}(q)
\\
\Delta_{3,2}(q)
\\
\Delta_{3,1}(q)
\end{array} \right) .
\en In other words, we have \be \notag
 {3\choose 1} \frac{\Gamma_{1}(q,3)\,
\Delta_{3,3}(q)}{3} = \Delta_{3,3}(q)\ , \en \be \notag
 {3\choose 2}
\frac{\Gamma_{2}(q,3)\, \Delta_{3,3}(q)}{3^2} + {2\choose 1}
\frac{\Gamma_{1}(q,3)\, \Delta_{3,2}(q)}{3}
 = \Delta_{3,2}(q)\ ,
\en \be {3\choose 3} \frac{\Gamma_{3}(q,3)\, \Delta_{3,3}(q)}{3^3} +
{2\choose 2} \frac{\Gamma_{2}(q,3)\,\Delta_{3,2}(q)}{3^2} +
{1\choose 1} \frac{\Gamma_{1}(q,3)\, \Delta_{3,1}(q)}{3}
 = \Delta_{3,1}(q)\ .
\notag \en
\end{ex}

 More generally, we have the following relation.

\begin{thm}
\label{thm S s}
 For any $1\leq k \leq n$ we have
 \be
\label{eq_colyrel} \Delta_{n,k}(q)  =  \sum_{m=0}^{n-k} {k+m\choose
1+m}
 \frac{\Gamma_{1+m}(q,n) \Delta_{n,k+m}(q)}{n^{1+m}}\ .
 \en
\end{thm}

\begin{proof}
Recall the notations $u_b$ and $D_b$ defined at the beginning of
Section \ref{subsec-derivsofseries}.
 Denote
\be
\pi_n =
  \prod_{b=1}^n\frac 1{1-u_b} ,
\qquad I_n^k(q) = \sum_{1\leq i_1<\dots < i_k\leq n}
\left(\prod_{j=1}^k
 \,(D_{i_j}+D_{i_j+1}+\dots + D_n)\right) \pi_n\ .
\en
 For example,
 \be
I^1_n(q) = (nD_n + (n-1)D_{n-1} + \dots + D_1)\pi_n . \notag
 \en
 Denote
 $\sigma_i({} c),\ i=1,\dots,n$, the elementary symmetric
functions in ${} c=(c_1,\dots,c_n)$ as before.

The left hand side of Formula (\ref{new n}) equals \be
 \sum_{k=0}^n\, k!(n-k)!\,\sigma_{n-k}({} c) I^k_n(q) .
\notag
 \en
 The right hand side equals
 \be
\frac {q^n}{n!}\, \sum_{k=0}^n\,
 k!(n-k)!\,\sigma_{n-k}({} c) \Delta_{n+1,n-k+1}(q)\ .
\notag
 \en
 By comparing coefficients of the symmetric polynomials in these two equations we get
\be \label{I s} I^k_n(q) = \frac {q^n}{n!}\,  \Delta_{n+1,n-k+1}(q)\
. \notag \en
  We have
\be \Delta_{n+1,n-k+1}(q) - \Delta_{n,n-k}(q) = (q-n)
 \Delta_{n,n-k+1}(q)
\notag \en
  and
  \be
I^k_n(q) = \sum_{j=0}^k{n-j \, \choose k-j}\,I^{j}_{n-1}(q) \,
\left((D_n)^{k-j}\frac 1{1-u_n} \right) . \notag \en
 Hence
 \be\label{Ink}
 I^k_n(q) - I^k_{n-1}(q)
\frac 1{1-u_n} = \sum_{j=0}^{k-1}{n-j\choose k-j}\,I^{j}_{n-1}(q) \,
\left((D_n)^{k-j}\frac 1{1-u_n} \right) . \en
 The left hand side of
 (\ref{Ink}) is
 \be
 \begin{aligned}
  I^k_n(q) - I^k_{n-1}(q)
\frac 1{1-u_n} &= I^k_n(q)
 - I^k_{n-1}(q)
\frac qn
\\
  &= \frac {q^n}{n!}\,  \Delta_{n+1,n-k+1}(q) - \frac
 {q^{n-1}}{(n-1)!}\,  \Delta_{n,n-k}(q) \frac qn\\
 &=\frac {q^n}{n!}\, ( \Delta_{n+1,n-k+1}(q) - \Delta_{n,n-k}(q))\\
  &= \frac
 {q^n(q-n)}{n!}\, \Delta_{n,n-k+1}(q) .
 \end{aligned}
 \notag
\en
  The right hand side of
 (\ref{Ink}) equals
 \be
 \begin{gathered}
 \sum_{j=0}^{k-1}{n-j\choose
 k-j}\,I^{j}_{n-1}(q) \, \left((D_n)^{k-j}\frac 1{1-u_n} \right) =
 \\
  \sum_{j=0}^{k-1} {n-j\choose k-j}\frac
{q^{n-1}}{(n-1)!}\,  \Delta_{n,n-j}(q)
\frac{q(q-n)}{n^{k-j+1}}\,\Gamma_{n-j}(q,n) .
 \end{gathered}
\notag
 \en
 These two equations
prove Theorem \ref{thm S s}.
\end{proof}

\medskip
\noindent
Since $\Delta_{n,n}(q) = 1$ and $\Gamma_1(q,n)=1$ for all $n$,
identity (\ref{eq_colyrel}) allows us
to express the polynomials $\Delta_{n,k}(q)$ in terms of polynomials
$\Gamma_k(q,n)$ and vice versa --- see the examples above.

\subsection{Identities involving Stirling numbers}
\label{subsec-stirlingIds}

 The conclusion of Theorem \ref{thm S s} can be written in the
 following form:

  \be
\label{eq_stirlingidentities}
  \Delta_{n,k}(q)=\frac{1}{n-k}\sum_{m=1}^{n-k}{m+k \choose m+1}\frac{1}{n^m}
  \Gamma_{m+1}(q,n)\Delta_{n,k+m}(q),
  \quad k=1, \dots,n.
  \en
  Comparing coefficients of $q^r$ of the two sides yields
  identities involving the Stirling numbers $s_1(n,k)$ and
  $s_2(n,k)$ of the first and second kinds,
  indexed with
$1 \leq k \leq n$ (and equalling 0
  outside this range).
  \begin{cor}
\label{cor-stirlingidentities}
  For
$1\leq k \leq n$,\ $0 \leq r \leq k$, we have
  \be
\label{coefs}
  \begin{gathered}
  {r+n-k-1 \choose r}
  s_1(n,r+n-k)=\\
  \sum_{m=0}^k {m+n-k \choose m+1}\sum_{i=0}^r
  {i+n-k+m-1 \choose i}\frac{(-1)^{m+r-i}}{n^{r-i+1}}(r-i+1)!\ \, \times \\
  s_2(m+1,r-i+1) s_1(n,i+n-k+m) .
  \end{gathered}
  \en
  \end{cor}

  For $r=0$, this says
  \be
  s_1(n,n-k)=\frac{n-k}{n}s_1(n,n-k)+\sum_{m=1}^k{m+n-k \choose m+1}
  \frac{(-1)^m}{n} s_1(n,n-k+m) ,
  \notag
\en and hence
  \be\label{star}
  s_1(n,n-k)=\frac{1}{k}\sum_{m=1}^k {m+n-k \choose m+1}(-1)^ms_1(n,n-k+m),
  \quad k=1, \dots, n.
  \en
  This contains the known formulas
  \be
  s_1(n,n-1)=-{n \choose 2}\quad\text{ and } \quad
  s_1(n,n-2)=\frac{1}{24}n(n-1)(n-2)(3n-1).
\notag
  \en

\end{document}